\documentclass{amsart}

\usepackage{amsmath} 
\usepackage{amssymb}
\usepackage{mathrsfs}

\newtheorem{theorem}{Theorem}[section] 
\newtheorem{claim}[theorem]{Claim}
 
\newtheorem{dc}[theorem]{Definition/Claim}

\theoremstyle{definition}
\newtheorem{definition}[theorem]{Definition}

\newtheorem{convention}[theorem]{Convention}

\theoremstyle{remark}
\newtheorem{remark}[theorem]{Remark}
\newtheorem{question}[theorem]{Question}

\newcommand{\rest}{{\restriction}}
 
\newcommand{\Hom}{{\rm Hom}}

\newcommand{\Rang}{{\rm Rang}}
\newcommand{\rk}{{\rm rk}}

\newcommand{\wilog}{{\rm without loss of generality}}

\newcommand{\then}{{\underline{then}}}
\newcommand{\when}{{\underline{when}}}
\newcommand{\Then}{{\underline{Then}}}

\newcommand{\mn}{{\medskip\noindent}}
\newcommand{\sn}{{\smallskip\noindent}}

\newcommand{\cH}{{\mathscr H}}

\newcommand{\cG}{{\mathscr G}}

\newcommand{\bbZ}{{\mathbb Z}}

\newcommand{\cP}{{\mathscr P}}

\newcommand{\cf}{{\rm cf}}

\newcount\skewfactor
\def\mathunderaccent#1#2 {\let\theaccent#1\skewfactor#2
\mathpalette\putaccentunder}
\def\putaccentunder#1#2{\oalign{$#1#2$\crcr\hidewidth
\vbox to.2ex{\hbox{$#1\skew\skewfactor\theaccent{}$}\vss}\hidewidth}}

\newenvironment{PROOF}[2][\proofname.]
   {\begin{proof}[#1]\renewcommand{\qedsymbol}{\textsquare$_{\rm #2}$}}
   {\end{proof}}

\begin{document}

\title {The First almost free Whitehead group}
\author {Saharon Shelah}
\address{Einstein Institute of Mathematics\\
Edmond J. Safra Campus, Givat Ram\\
The Hebrew University of Jerusalem\\
Jerusalem, 91904, Israel\\
 and \\
 Department of Mathematics\\
 Hill Center - Busch Campus \\ 
 Rutgers, The State University of New Jersey \\
 110 Frelinghuysen Road \\
 Piscataway, NJ 08854-8019 USA}
\email{shelah@math.huji.ac.il}
\urladdr{http://shelah.logic.at}
\thanks{Research supported by German-Israeli Foundation for Scientific
Research/Development Grant No. I-706-54.6/2001. Publication 914.
\newline
I would like to thank Alice Leonhardt for the beautiful typing.}



\date{June 7, 2011}

\begin{abstract}
Assume G.C.H. and $\kappa$ is the first uncountable
cardinal such that there is a non-free $\kappa$-free abelian Whitehead group of
cardinality $\kappa$.  We prove that if all $\kappa$-free Abelian
group of cardinality $\kappa$ are Whitehead then $\kappa$ is necessarily an
inaccessible cardinal.
\end{abstract}

\maketitle
\numberwithin{equation}{section}
\setcounter{section}{-1}

\section{Introduction} 

For the non-specialist reader, note that we deal exclusively with
abelian groups, we call such $G$ Whitehead \when \,: if $\bbZ
\subseteq H$ and $H/\bbZ \cong G$ then $\bbZ$ is a direct summand
of $H$.
Why this property is worthwhile and  generally on the subject see the
book of Eklof-Mekler \cite{EM02}.  Recall $G$ is a free abelian group
if it is the direct sum of copies of $(\bbZ,+)$.
All free abelian groups are
Whitehead and every Whitehead group is ``somewhat free" ($\aleph_1$-free, 
see below).  Possibly (i.e. consistently with ZFC) every Whitehead 
group is free and possibly not (in fact, it seems that almost 
any behaviour is possible).

Recall that G.C.H., the generalized continuum hypothesis, says that
$2^\lambda = \lambda^+$ for every infinite $\lambda$ and ``$\lambda$
is an inaccessible cardinal" means that: it is strong limit (i.e. $\mu <
\lambda \Rightarrow 2^\mu < \lambda$) and regular, i.e. $\lambda >
\sum\limits_{i < \kappa} \lambda_i$: when $\kappa < \lambda$ and $i <
\kappa \Rightarrow \lambda_i < \lambda$.

It is well known that inaccessible cardinals are large, e.g. it is not
provable in ZFC, the usual axioms of set theory that such cardinals exists.

Also being ``$\kappa$-free of cardinality $\kappa$" is a central
 notion where we say an abelian group 
$G$ is $\kappa$-free when the pure closure inside $G$ of any 
subgroup generated by $<
 \kappa$ elements is free.  This explains the interest in the 
result, which restricts the behaviour.
That is, assume $\kappa$ is the first cardinal $\lambda$ such that
there is a non-free $\lambda$-free abelian Whitehead group of cardinality
$\lambda$.  The conclusion is, assuming GCH and $\kappa$ not
 strongly inaccessible, that not all such abelian groups are Whitehead. 

Set theorists may recall that in 
\cite[\S1]{Sh:667} it is proved that if $\mu$ is strong limit
singular, $\lambda = \mu^+ = 2^\mu$ and $S \subseteq \{\delta <
\lambda:\text{cf}(\delta) = \text{cf}(\mu)\}$ is stationary, \then \,
though $\diamondsuit_S$ may (= consistently) fail, still we can prove  
a relative of $\diamondsuit_S$ sufficient for constructing 
abelian groups $G$ of cardinality $\lambda$ related to 
satisfying Hom$(G,\bbZ) \ne \{0\}$.  We prove here in \ref{5n.2} somewhat
more and use it in the proof of the Theorem \ref{5n.1}.  This
claim is complementary to \cite{Sh:587} which shows that consistently
there is such regular (in fact strongly inaccessible) $\lambda$.
\bigskip

\begin{convention}
\label{5n.0A}
Saying group we mean abelian group.
\end{convention}

\noindent
In \cite{Sh:587} we answer
\begin{question}
\label{5n.0D}
(G\"obel) Assuming GCH can there be a regular $\kappa$ such that:
\mn
\begin{enumerate}
\item[$\boxdot_\kappa$]  $(a) \quad \kappa = \text{ cf}(\kappa) >
  \aleph_0$ and there are $\kappa$-free not free groups of cardinality $\kappa$
\sn
\item[${{}}$]  $(b) \quad$ every $\kappa$-free group of cardinality
$\kappa$ is a Whitehead group.
\end{enumerate}
\mn
Moreover by \cite{Sh:587} it is consistent that: G.C.H. + for some
strongly inaccessible $\kappa$ we have $\circledast_\kappa$ where the
statement $\circledast_\kappa$ is defined by:
\mn
\begin{enumerate}
\item[$\circledast_\kappa$]  $(a) \quad \kappa$ is regular
uncountable and there are $\kappa$-free non-free (abelian) groups

\hskip25pt   of cardinality $\kappa$
\sn
\item[${{}}$]   $(b) \quad$ every $\kappa$-free (abelian) group of
cardinality $\kappa$ is a Whitehead group
\sn
\item[${{}}$]  $(c) \quad$ every Whitehead group of cardinality $<
\kappa$ is free.
\end{enumerate}
\mn
[For clause (c), the following sufficient condition is used there: for every 
regular uncountable $\lambda <
\kappa$ we have $\diamondsuit^*_\lambda$.  The proof starts with
$\kappa$ weakly compact and adds no new sequence of length $< \kappa$,
so e.g. starting there with $\bold L$ this holds.]
\end{question}

\noindent
A natural question is whether
\begin{question}
\label{5n.0H}
Assume G.C.H., is it consistent that there is an accessible $\kappa$
such that $\circledast_\kappa$ holds?  In other words there
is cardinal $\kappa$ satisfying the following and there first such $\kappa$
is accessible:
\mn
\begin{enumerate}
\item[$(*)$]  there is a $\kappa$-free 
non-free Whitehead group of cardinality $\kappa$.
\end{enumerate}
\end{question}

\noindent
Now Theorem \ref{5n.1} says that no.
But another natural question is:

\begin{question}
\label{5n.0K}
Assume G.C.H.  Can there be a $\kappa$
such that $\boxdot_\kappa$ but is the first such $\kappa$ accessible?  

This is problem F4 of \cite{EM02} and it remains open.
\end{question}

\begin{definition}
\label{0z.11}
1) An abelian group $G$ is free if $G$ is the direct sum of copies of $\bbZ$.

\noindent
2) For a cardinality $\kappa \ge \aleph_1$, an 
abelian group $G$ is $\kappa$-free
\when \, every subgroup of $G$ of cardinality $< \kappa$ is free.

\noindent
3) $\bar G = \langle G_\alpha:\alpha < \lambda\rangle$ is a filtration
of the abelian group $G$ of cardinality $\lambda$ if $G_\alpha$ is
a subgroup of $G$ of cardinality $< \lambda$, increasing continuous
with $\alpha$ and $G = \cup\{G_\alpha:\alpha < \lambda\}$.

\noindent
4) $G$, an abelian group of cardinality $\kappa$, is strongly
$\kappa$-free \when \, it is $\kappa$-free and for every subgroup
$G'$ of $G$ of cardinality $< \kappa$ there is a subgroup $G''$ of
cardinality $< \kappa$ such that $G' \subseteq G''$ and $G/G''$ is
 $\kappa$-free. 
\end{definition}

\noindent
We thank the referees for many helpful comments, in particular for
making the work more self-cointained.
\newpage

\section {The First almost free non-free Whitehead} 

\begin{theorem}
\label{5n.1}
(G.C.H.)  Let $\kappa$ be the first
$\lambda > \aleph_0$ such that there is a $\lambda$-free abelian 
Whitehead group, not free, of cardinality $\lambda$ (and we assume 
there is such $\lambda$).  If $\kappa$ is not (strongly) inaccessible, 
\then \, there is a non-Whitehead group $G$ of cardinality $\kappa$
which is $\kappa$-free (and necessarily non-free). 
\end{theorem}

\begin{PROOF}{\ref{5n.1}}
\medskip

\noindent
\underline{Stage A}:  Let $G$ exemplify the choice of $\kappa$.

Necessarily $\kappa$ is regular, by the singular compactness, see
e.g. \cite[Ch.IV,3.5]{EM02}.  Let $\bar G = \langle
G_\alpha:\alpha < \kappa\rangle$ be a filtration of $G$ so as $G$ is
$\kappa$-free necessarily every $G_\alpha$ is free.  Without loss
of generality each $G_\alpha$ is a pure subgroup of $G$.  Let $S :=
\Gamma(\bar G) = \{\alpha < \kappa:\alpha$ is a limit ordinal 
and $G/G_\alpha$ is not $\kappa$-free$\}$.  Now recall $G$ is not free hence 
$\Gamma(\bar G)$ is stationary.
\medskip

\noindent
\underline{Stage B}:  $G$ is strongly $\kappa$-free.

Toward contradiction assume that not.  Without loss of generality 
$\kappa$ is a successor cardinal (why? by the theorem's assumption and
$\kappa$ being regular but also as for $\kappa$ a limit cardinal, 
$\kappa$-free implies strongly $\kappa$-free).

By our present assumption toward contradiction, for some $\alpha <
\kappa$ for every $\beta \in [\alpha,\kappa)$, the abelian group,
$G/G_\beta$ is not $\kappa$-free.   Hence \wilog \, $\alpha < \lambda
\Rightarrow G_{\alpha +1}/G_\alpha$ is not free
and if $G_{\alpha +1}/G_\alpha$ is
uncountable then for some $\kappa_\alpha$, it is $\kappa_\alpha$-free
not free.  By the theorem assumption (and as countable Whitehead
groups are free), $\alpha < \kappa \Rightarrow G_{\alpha +1}/G_\alpha$ is
not a Whitehead group, i.e. $\Gamma(\bar G) = \kappa$.  But $\kappa$ is a
successor cardinal, so let $\kappa = \mu^+$; also $2^\mu < 2^\kappa$
as we are assuming GCH so the weak diamond holds for $\kappa$ (see
\cite{DvSh:65}) so by the previous sentence and
e.g. \cite[1.10,pg.369]{EM02},
we know that $G$ is not a Whitehead group, contradiction to the
assumption on $G$.
\medskip

\noindent
\underline{Stage C}:   Hence \wilog
\mn
\begin{enumerate}
\item[$(*)_1$]   if $\alpha$ is a non-limit ordinal then $G/G_\alpha$ is
$\kappa$-free
\end{enumerate}
and obviously \wilog \,
\mn
\begin{enumerate}
\item[$(*)_2$]   $(a) \quad$ if $\alpha \in S$ then 
$G_{\alpha +1}/G_\alpha$ is not free
\sn
\item[${{}}$]   $(b) \quad$ if $\alpha < \kappa$ then $G_\alpha$
is a pure subgroup of $G$ and is free 

\hskip25pt hence $G/G_\alpha$ is torsion free.
\end{enumerate}
\mn
Also we can choose $\bar H = \langle H_\alpha:\alpha \in S\rangle$ such
that
\mn
\begin{enumerate}
\item[$(*)_3$]   $(a) \quad H_\alpha$ is a subgroup of $G$
\sn
\item[${{}}$]  $(b) \quad H_\alpha/(H_\alpha \cap G_\alpha)$ is not free
\sn
\item[${{}}$]  $(c) \quad$ under (a) + (b) the rank of $H_\alpha/(H_\alpha
\cap G_\alpha)$ is minimal call it $\theta_\alpha$ hence

\hskip25pt $\theta_\alpha$ is $< \aleph_0$ or is regular uncountable $< \kappa$
\sn
\item[${{}}$]  $(d) \quad$ the cardinality of $H_\alpha$ is $\le
  \theta_\alpha + \aleph_0$, in fact equality holds.
\end{enumerate}
\mn
Note that $|H_\alpha|$ may be $< |G_\alpha|$ so it is unreasonable to
ask for $G_\alpha \subseteq H_\alpha$. 
\mn
\begin{enumerate} 
\item[$(*)_4$]   without loss of generality
\begin{enumerate}
\item[$(a)$]   $H_\alpha \subseteq G_{\alpha +1}$
\sn
\item[$(b)$]   $G_\alpha/(H_\alpha \cap G_\alpha)$ is free
\sn
\item[$(c)$]  $G_\alpha + H_\alpha$ is a pure subgroup of $G_{\alpha +1}$.
\end{enumerate}
\end{enumerate}
\mn
[Why?  For clause (a) we can restrict $\bar G$ to a club.  
For clause (c) let $H''_\alpha$ be the pure closure of $H_\alpha + 
G_\alpha$ inside $G_{\alpha +1}$ so $H''_\alpha/G_\alpha,(H_\alpha +
G_\alpha)/G_\alpha,H_\alpha/(H_\alpha \cap G_\alpha)$ has
the same rank, which is $\theta_\alpha$.  As $G_{\alpha +1}/G_\alpha$
is torsion free, also $H''_\alpha/G_\alpha$ is torsion free.  Hence
there is $H'_\alpha \subseteq H''_\alpha$ of cardinality $\le
\theta_\alpha + \aleph_0$ such that $G_\alpha + H'_\alpha$ is a pure
subgroup of $H''_\alpha$ hence of $G_{\alpha +1}$ and
$H'_\alpha + G_\alpha = H''_\alpha$ and e.g. the rank of
$H'_\alpha/(H'_\alpha \cap G_\alpha)$ is the same as the rank of
$H''_\alpha/G_\alpha$ which is $\theta_\alpha$, so replacing
$H_\alpha$ by $H'_\alpha$ also clause (c) holds.

For clause (b) note that $G_\alpha$ is free; so
there is $G'_\alpha \subseteq G_\alpha$ of cardinality $\theta_\alpha
+ \aleph_0$ such that $G_\alpha/G'_\alpha$ is free and $H_\alpha \cap G_\alpha
\subseteq G'_\alpha$ and replace $H_\alpha$ by $H^*_\alpha = 
H_\alpha + G'_\alpha$ noting that $H^*_\alpha + G_\alpha = H_\alpha +
G_\alpha$.] 

By the hypothesis on $\kappa$, if $K$ is a $\lambda$-free not free
group of cardinality $\lambda$ and $\aleph_0 < \lambda < \kappa$ then
$K$ is not a Whitehead group so clearly (recalling that Whitehead
groups which are countable are free):
\mn
\begin{enumerate}
\item[$(*)_5$]   $H_\alpha/(G_\alpha \cap H_\alpha)$ is not Whitehead
  for $\alpha \in S$.
\end{enumerate}
\mn
Hence by \cite{Sh:44} or see \cite[Ch.VI,1.13]{EM02} we know that
\mn
\begin{enumerate}
\item[$(*)_6$]   $\neg \diamondsuit_S$.
\end{enumerate}
\mn
Recall $\kappa$ is regular uncountable so toward contradiction assume
$\kappa = \mu^+$.  Let $\sigma = \text{ cf}(\mu)$.  But,
\cite{Sh:108} or \cite[Ch.VI,1.13]{EM02}, GCH implies that $\diamondsuit_{S'}$ 
for every stationary $S' \subseteq \kappa \backslash S^\kappa_\sigma$
where $S^\kappa_\sigma := \{\delta < \kappa:\text{cf}(\delta) = \sigma\}$.

Hence we know that
\mn
\begin{enumerate}
\item[$(*)_7$]   for some club $E$ of $\kappa$ we have $S \cap E \subseteq
S^\kappa_\sigma$ so \wilog \, $S \subseteq S^\kappa_\sigma$, i.e.,
$\delta \in S \Rightarrow \text{ cf}(\delta) = \sigma$.
\end{enumerate}
\mn
Also \wilog \, (from ``strongly $\kappa$-free")
\mn
\begin{enumerate}
\item[$(*)_8$]   $\delta \in S \Rightarrow \theta_\delta + \aleph_0 
\ge \sigma$.
\end{enumerate}
\mn
[Why?  For completeness we elaborate.
Let $S^1 = \{\delta \in S:\theta_\delta + \aleph_0 < \sigma\}$;
  first assume $S^1$ is stationary.  For each $\delta \in S^1$
  necessarily $H_\delta \cap G_\delta$ has cardinality $\le
  \theta_\delta + \aleph_0 < \sigma = \cf(\delta)$ hence for some
  $\alpha_\delta < \delta,H_\delta \cap G_\delta \subseteq
  G_{\alpha_\delta}$.  By Fodor lemma for some $\alpha(*) < \kappa$
  the set $S^2 := \{\delta \in S^1:\alpha_\delta = \alpha(*)\}$ is
  stationary.  As $\sigma = \cf(\mu),\kappa = \mu^+$, clearly $\mu^{<
    \sigma} = \mu$ recalling GCH hence $\{H_\delta \cap
  G_{\alpha(*)}:\delta \in S^2\}$ has cardinality $\le \mu$ hence $S^3
  = \{\delta \in S^2:H_\delta \cap G_{\alpha(*)} = H_*,\theta_\delta =
  \theta_*\}$ is a stationary subset of $\kappa$ for some $\theta_*$
  and $H_* \subseteq H_{\alpha(*)}$.  Let $\alpha_\varepsilon$ be the
  $\varepsilon$-th member of $S^3$ and $H^\varepsilon = H_* +
  \Sigma\{H_{\alpha_\zeta}:\zeta < \varepsilon\}$ for $\varepsilon \le
  (\theta_* + \aleph_0)^+$.  Clearly for $\varepsilon = (\theta_* +
  \aleph_0)^+,H^\varepsilon$ is not free and has cardinality $\le
  \sigma = \cf(\mu) < \lambda$, contradiction to ``$G$ is
  $\kappa$-free".  So necessarily $S^1$ is not stationary, and as we
  can restrict $\bar G$ to a club, \wilog \, $S^1 = \emptyset$ so
  $(*)_8$ holds indeed.]
\medskip

\noindent
\underline{Stage D}:  $\alpha \in S \Rightarrow \theta_\alpha < \mu$.

\noindent
Why?  Otherwise for some $\alpha \in S,\theta_\alpha = \mu$ so
 $\theta_\alpha$ is infinite hence 
$\theta_\alpha$ is regular and by $(*)_3$ uncountable and there is a
 $\theta_\alpha$-free not free group of cardinality $\theta_\alpha$ hence
 $\mu$ is regular and there is a $\mu$-free not free
abelian group of cardinality $\mu < \kappa$, i.e. $H_\alpha/(G_\alpha
\cap H_\alpha)$, hence this group is not Whitehead.  So there
is a sequence $\langle H^*_\varepsilon:\varepsilon \le \mu +1\rangle$
purely increasing continuous sequence of free groups
such that $\varepsilon < \mu$ implies
$H^*_{\mu+1}/H^*_\varepsilon$ is free but $H^*_{\mu
+1}/H^*_\mu$ is not free and is isomorphic to $H_\alpha/(G_\alpha \cap
H_\alpha)$. 

[Why?  Let $\langle y_\varepsilon:\varepsilon < \mu\rangle$ be a free
  basis of $G_\alpha$, so $G_\alpha = \oplus\{\bbZ
  y_\varepsilon:\varepsilon < \mu\}$.
As $\{y_\varepsilon:\varepsilon < \mu\} \subseteq G_\alpha$ and $\cf(\alpha) =
\mu$, there is an increasing continuous sequence $\langle
\gamma_\varepsilon:\varepsilon < \mu\rangle$ of ordinals $< \alpha$
with limit $\alpha$ such that $y_\varepsilon \in
G_{\gamma_{\varepsilon +1}}$.

Let $H^*_\varepsilon = \oplus\{\bbZ y_\zeta:\zeta < \varepsilon\}$ for
$\varepsilon < \mu,H^*_\mu = G_\alpha,H^*_{\mu +1} = G_\alpha + H_\alpha$,
they are as required.  E.g. why is $H^*_{\mu +1}/H_\varepsilon$ 
free? as $G_{\alpha +1}/G_{\alpha_\varepsilon +1}$ is free by $(*)_1$,
also $G_\alpha/H^*_\varepsilon$ is free as $\{y_\zeta +
H^*_\varepsilon:\zeta \in [\varepsilon,\mu)\}$ is a free basis but
  $H^*_\varepsilon \subseteq G_{\alpha_{\varepsilon +1}} \subseteq
  G_\alpha$ hence $G_{\alpha_\varepsilon +1}/H^*_\varepsilon$ is free
  so together $G_{\alpha +1}/H^*_\varepsilon$ is free which implies
  that its subgroup $H^*_{\mu +1}/H_\varepsilon$ is free as promised.]

We can find $H_*$ and $\langle H_\eta,h_\eta:\eta \in {}^{\mu \ge}2\rangle$
such that:
\mn
\begin{enumerate}
\item[$(a)$]   $H_* = \oplus\{\bbZ x_t:t \in I\}$ and $|I|=\mu$
\sn
\item[$(b)$]   $I$ is the disjoint union of $\{I_\eta:\eta \in{}^{\mu >}2\}$
\sn
\item[$(c)$]  $|I_\eta| = \rk(H^*_{\ell g(\eta)+1}/H^*_{\ell
g(\eta)})$ for $\eta \in {}^{\mu >} 2$
\sn
\item[$(d)$]   $H_\eta = \oplus\{\bbZ x_t:t \in \cup\{I_{\eta
\rest \varepsilon}:\varepsilon < \ell g(\eta)\}\} \subseteq H_*$ for
$\eta \in {}^{\mu >} 2$
\sn
\item[$(e)$]  $h_\eta$ is an isomorphism from $H_{\ell g(\eta)}$
onto $H_\eta$
\sn
\item[$(f)$]  $h_{\eta \rest \varepsilon} \subseteq h_\eta$ if
$\varepsilon < \ell g(\eta),\eta \in {}^{\mu >} 2$.
\end{enumerate}
\mn
Now we can find $\langle H^+_\eta,h^+_\eta:\eta \in {}^\mu 2\rangle$
such that
\mn
\begin{enumerate}
\item[$(g)$]   $H_\eta \subseteq H^+_\eta$ 
\sn
\item[$(h)$]  $h^+_\eta$ is an isomorphism from $H^*_{\mu + 1}$
onto $H^+_\eta$ extending $h_\eta$.
\end{enumerate}
\mn
Without loss of generality
\mn
\begin{enumerate}
\item[$(i)$]   $H^+_\eta \cap H_* = H_\eta$
\end{enumerate}
\mn
so there is
\mn
\begin{enumerate}
\item[$(j)$]  $H^*_\eta$ extends $H^+_\eta$ and $H_*$ such that $H^+_\eta
\cup H_*$ generates $H^*_\eta$ and $H^*_\eta/H_*$ is isomorphic to
$H^*_\eta/H_\eta$. 
\end{enumerate}
\mn
Lastly, \wilog \,
\mn
\begin{enumerate}
\item[$(k)$]   the sets $\langle H^*_\eta \backslash H_*:\eta \in
{}^\mu 2\rangle$ are pairwise disjoint, 
\end{enumerate}
\mn
so there is an abelian group $H$ such that
\mn
\begin{enumerate}
\item[$(l)$]   $H^*_\eta \subseteq H$ for $\eta \in {}^\mu 2$ and
\sn
\item[$(m)$]   $H/H_* = \oplus\{H^*_\eta/H_*:\eta \in {}^\mu 2\}$
\sn
\item[$(n)$]   $H$ has cardinality $2^\mu = \kappa$.
\end{enumerate}
\mn
Next we note
\mn
\begin{enumerate}
\item[$(o)$]  $H$ is $\kappa$-free.
\end{enumerate}
\mn
[Why?  Note that $\cup\{H^+_\eta:\eta \in {}^\mu 2\}$ include $\{x_t:t
  \in I\}$ hence the subgroup it generates includes $H$.  
Let $H' \subseteq H$ be a sub-group
of cardinality $< \kappa$ so $\le \mu$, hence there are $\eta_i
\in {}^\mu 2$ for $i < \mu$ such that $H' \subseteq
\Sigma\{H^+_{\eta_i}:i < \mu\}$ and $j < i \Rightarrow \eta_i \ne
\eta_j$.  We can choose $\zeta(i) < \mu$ by induction on $i < \mu$
such that $\langle \{\eta_i \rest \varepsilon:\varepsilon \in
[\zeta(i),\mu)\}:i < \mu\}$ is a sequence of pairwise disjoint sets.

For each $i$ let $H^{**}_i \subseteq H^+_{\eta_i}$
be such that $H^+_{\eta_i} = H^{**}_i \oplus H_{\eta_i \rest
  \zeta(i)}$.  
Let us define $H'_i$ for $i \le \mu$ by: $H'_0 \subseteq H$ is
generated by $\{x_t:t \in I$ but $t \notin
\cup\{I_{\eta_i \rest \varepsilon}:i < \mu$ and $\varepsilon \in
[\zeta(i),\mu)\}\}$ and $H'_i$ is $\oplus\{H^{**}_j:j<i\} \oplus
  H'_0$.  Clearly $H' \subseteq H'_\mu \subseteq H,\langle H'_i:i \le
  \mu\rangle$ is increasing continuous, $H'_0$ is free and
  $H'_{i+1}/H'_i \cong H^{**}_i$ is free.  It follows that $H'_\mu$ is
  free hence $H' \subseteq H'_\mu$ is free.  So clause (o) holds indeed.]

Let $\langle \eta_\alpha:\alpha < \kappa\rangle$ list ${}^\mu 2$, let
$H^{**}_\alpha \subseteq H$ be $\Sigma\{H^+_{\eta_\beta}:\beta <
\alpha\} + H_*,\bar H^{**} = \langle H^{**}_\alpha:\alpha <
\kappa\rangle$ is a filtration of $H$, and $\Gamma(\bar
H^{**})=\kappa$ and $H^{**}_{\alpha +1}/H^{**}_{\alpha +1}$ is
isomorphic to $H^*_{\mu +1}/H^*_\mu$ so is not free and is not
Whitehead.  Hence, see \cite[Ch.XII,1.10,pg.369]{EM02}, $H$ is not 
Whitehead.  So $H$ is $\kappa$-free (see clause (o) of
cardinality $\kappa$ (see clause (o)  and not Whitehead, and so
not free), the desired conclusion of the theorem.  So indeed \wilog \,
the stage desired conclusion, $\alpha \in S \Rightarrow \theta_\alpha
< \mu$ holds.]
\medskip

\noindent
\underline{Stage E}:  $\mu$ is singular.

Why?  Because by earlier stages cf$(\mu) = \sigma$ and 
$\alpha \in S$ implies $\text{cf}(\alpha) \le \theta_\alpha < \mu$
and $\alpha \in S \Rightarrow \sigma = \text{ cf}(\alpha)$,
so necessarily $\mu$ is singular.
\medskip

\noindent
\underline{Stage F}:  Let $\sigma = \text{ cf}(\mu)$ so $\sigma$ is regular
$< \mu$; also choose $\theta = \text{ cf}(\theta)
< \mu$ such that $S^*_1 := \{\delta \in S:\theta_\delta = \theta\}$ is
stationary.  Note that $\mu = \mu^{< \sigma}$ as G.C.H. holds
recalling $\sigma = \cf(\delta)$ and $\delta \in S \Rightarrow
\cf(\delta) = \sigma$ by $(*)_7$.  We shall now use \cite[\S3]{Sh:521}
and its notation, see \ref{1n.4}, \ref{1n.6}, \ref{5n.2} below.

By the Theorem \ref{1n.4} below we can find a $\kappa$-witness
$\bold x$, so it consists of
$n \ge 1,\bold S,\langle B_\eta:\eta \in
\bold S_c\rangle,\langle s^\ell_\eta:\eta \in \bold S_f,\ell < n\rangle$ as
there, i.e. as in \cite[3.6]{Sh:521} with $(\lambda,\kappa^+,S)$ there standing
for $(\kappa,\aleph_1,\bold S)$ here, such that
\mn
\begin{enumerate}
\item[$(*)$]   $\langle \alpha \rangle \in \bold S \Leftrightarrow \alpha
\in S^*_1$, this means $W(<>,\bold S) = S^*_1$.
\end{enumerate}
\mn
Clearly continuing to use the notation there, $\alpha \in S^*_1 \Rightarrow
\lambda(\langle \alpha\rangle,\bold S) \le \mu$ hence being
regular it is $< \mu$, and so \wilog \, constant (in fact it is
$\theta$ by the proof).  Now we apply the
claim \ref{5n.2} below with $\lambda$ there standing for $\kappa$
here.  

Why clause (d) of the assumption of \ref{5n.2} holds?  If $\theta >
\aleph_0$, by \cite[1.2]{Sh:521}, the group $G_{\bold x(<\alpha>)}$,
derived from $\bold x(\langle \alpha \rangle)$ is a $\lambda(<
\infty,\bold S_{\bold x})$-free non-free abelian group of cardinality
$\lambda(\langle \alpha \rangle,\bold S_{\bold x}) = \theta$, so is
not Whitehead by the theorem assumption on $\kappa$; note that
$G_{\bold x(< \alpha>)}$ was derived in some way from
$H_\alpha/(H_\alpha \cap G_\alpha)$, but it is not necessarily equal
to it.  If $\theta \le \aleph_0$ then $G_{\bold x(<\alpha>)}$ is a
non-free (abelian) countable group hence is not Whitehead.

So the assumption of \ref{5n.2} says that 
there is a strongly $\kappa$-free abelian group $G$ of cardinality
$\kappa$ by the theorem's assumption which is
not Whitehead, so we are done. 
\end{PROOF}

\noindent
Recall (by \cite[\S5]{Sh:161} and see \cite[\S3]{Sh:521}; or see \cite{EM02}).
\begin{theorem}
\label{1n.4}
For any $\lambda > \aleph_0$ the following conditions are equivalent:
\mn
\begin{enumerate}
\item[$(a)$]   there is a $\lambda$-free not free abelian group
\sn
\item[$(b)$]   {\rm PT}$(\lambda,\aleph_1)$ which means that: there
is a family ${\cP}$ of countable sets of cardinality $\lambda$ with
no transversal (i.e. a one-to-one choice function) but any subfamily of
cardinality $< \lambda$ has a transversal
\sn
\item[$(c)$]  there is a $\lambda$-witness $\bold x$ as in
\cite[3.6,3.7]{Sh:521}; so $\bold x$ consists of
\sn
\begin{enumerate}
\item[$(\alpha)$]   a natural number $n$
\sn
\item[$\beta)$]   a so-called $\lambda$-set $\bold S 
\subseteq \{\eta \in {}^{n \ge}\lambda:\eta$
decreasing$\}$ closed under initial segments, (see \cite[3.1]{Sh:521})
\sn
\item[$(\gamma)$]  disjoint $\lambda$-system $\bar B = 
\langle B_\eta:\eta \in \bold S_c\rangle$, see \cite[3.4]{Sh:521}
\sn
\item[$(\delta)$]   $\bar s = \langle s^\ell_\eta:\eta
\in \bold S_f,\ell < n\rangle$, etc.
\end{enumerate}
\item[$(d)$]  Without loss of generality $\cup\{B_\eta:\eta \in \bold S_c\} =
\lambda$, so $\bold S = \bold S_{\bold x}$, etc.
\sn
\item[$(e)$]    Let $\langle
a^\ell_{\eta,m}:m < \omega\rangle$ list $s^{\ell,\bold x}_\eta$ with
no repetition for $\eta \in \bold S^{\bold x}_f,\ell < n$ and
$a^{\ell(1)}_{\eta(1),m(1)} = a^{\ell(2)}_{\eta(2),m(2)}$ implies
$\ell(1) = \ell(2),m(1) = m(2)$ and $m<m(1) \Rightarrow
a^{\ell(1)}_{\eta(\ell),m} = a^{\ell(2)}_{\eta(2),m}$, so \wilog \,
$a^\ell_{\eta,m} < a^\ell_{\eta,m+1}$ for every relevant $\eta,\ell,m$;
also $\delta = \sup \cup \{s^{0,\bold x}_\eta:\langle \delta \rangle
\trianglelefteq \eta \in \bold S_f\}$ when $\langle \delta \rangle \in
\bold S_f$ and $\eta \mapsto \lambda(\eta,\bold S_{\bold x}),
\eta \mapsto W(\eta,S_{\bold x})$ are well defined.
\end{enumerate}
\end{theorem}

\noindent
Used but not named in \cite{Sh:521}:
\begin{dc}
\label{1n.6}
1) For a $\lambda$-witness
$\bold x$ and $\nu \in \bold S_{\bold x}$ we define $\bold x_\nu =
\bold x(\nu)$ by $n_{\bold x(\nu)} = n_{\bold x} - \ell
g(\nu),\bold S_{\bold x(\nu)} = \{\eta:\nu \char 94 \eta \in
\bold S_{\bold x}\},B^{\bold x(\nu)}_\eta = B^{\bold x}_{\eta
\char 94 \nu},s^{\ell,\bold x(\nu)}_\eta = s^{\ell g(\nu) +
\ell,\bold x}_{\nu \char 94 \eta}$.

\noindent
2) For a $\lambda$-witness $\bold x$ let $G_{\bold x}$ be the abelian
group $G_{\{<\alpha>:\alpha \in W(<>,S_{\bold x})\}}$ defined inside
the proof of \cite[1.2]{Sh:521}.
\end{dc}

\begin{remark}
\label{a8}
We may use the following.

\noindent
1) For a $\lambda$-witness $\bold x$ let 

\begin{equation*}
\begin{array}{clcr}
K_{\bold x} = \{I \subseteq 
\bold S_{\bold x}:\, I &\text{ is a set of pairwise }
\triangleleft\text{-incomparable sequences such that} \\
  &\{\beta:\eta \char 94 \langle \beta \rangle \in I\} \text{ is an
  initial segment of} \\
  &W(\eta,S) \text{ for any } \eta \in \bold S\}.
\end{array}
\end{equation*}

\mn
2) For $\bold x$ a $\lambda$-witness and $I \in K_{\bold x}$ let
   $Y[I]$ and $G_I$ be defined as in the proof of \cite[1.2]{Sh:521}
   before Fact A.  We may write $\eta$ instead of $I = \{\eta\}$.
\end{remark}

\begin{claim}
\label{5n.2}
Assume
\mn
\begin{enumerate}
\item[$(a)$]  $\mu$ is strong limit singular, $\lambda = \mu^+ =
2^\mu$ and $\sigma = \,\text{\rm cf}(\mu)$
\sn
\item[$(b)$]  $S \subseteq \{\delta < \lambda:\text{\rm cf}(\delta)
= \sigma\}$ is stationary
\sn
\item[$(c)$]  $\bold x$ is a $\lambda$-witness see \ref{1n.4} with
$W(\langle \rangle,\bold S) \subseteq S$
\sn
\item[$(d)$]  for each $\alpha \in W(\langle \rangle,\bold S)$, the abelian
group $G_{\bold x(<\alpha>)}$ is not Whitehead
(where $G_{\bold x(<\alpha>)}$ is defined as
inside the proof of \cite[1.2]{Sh:521}).
\end{enumerate}
\mn
\Then \,

\noindent
1) There is a strongly $\lambda$-free abelian group $G$ of cardinality
 $\lambda$ which is not Whitehead, in fact $\Gamma(G) \subseteq S$.

\noindent
2) There is a strongly $\lambda$-free abelian group $G^*$ of cardinality
$\lambda$ satisfying {\rm HOM}$(G^*,\bbZ) = \{0\}$, in fact
$\Gamma(G^*) \subseteq S$ (in fact the same abelian group can serve).
\end{claim}

\begin{remark}
1) We rely on \cite[\S3]{Sh:521}.

\noindent
2) So in clause (d), $G_{\bold x(<\alpha>)}$ is the abelian group 
defined from $\bold x(\langle \alpha \rangle)$.

\noindent
3) If you do not like clause (d) of \ref{5n.2}, replace it by
``$\lambda$ is as in \ref{5n.1}".
\end{remark}

\begin{PROOF}{\ref{5n.2}}
1) Let $\bold S = \bold S_{\bold x}$, etc.  Without loss of generality
\mn
\begin{enumerate}
\item[$(*)_0$]   $\cup\{B^{\bold x}_\eta:\eta \in \bold S_c\} = \lambda$.
\end{enumerate}
\mn
Let ${\cH}(\lambda) = \bigcup\limits_{\alpha < \lambda}
M_\alpha$ where $M_\alpha \prec ({\cH}(\lambda),\in)$ has
cardinality $\mu$, is increasing continuous with $\alpha$ such that
$\mu +1 \subseteq M_0$ and $\langle
M_\beta:\beta \le \alpha\rangle \in M_{\alpha +1}$.

Let $S_0 = \{\delta \in W(\langle \rangle,\bold S):M_\delta \cap \lambda =
\delta\}$, as $W(\langle \rangle,\bold S)$ is stationary and $\{\delta
< \lambda:M_\delta \cap \lambda = \delta\}$ is a club of $\lambda$; clearly
 also $S_0$ is a stationary subset of $\lambda$; now for each $\delta \in S_0$ 
as $\mu$ is strong limit (and so $\mu = \mu^{< \sigma}$) clearly
\mn
\begin{enumerate}
\item[$(*)_1$]   if $Y \subseteq M_\delta,|Y| < \mu$ and cf$(\delta)
= \sigma$ \then \, there is
an increasing continuous sequence $\langle X_i:i < \sigma\rangle$ of
subsets of $Y$ with union $Y$ such that $M_\delta \supseteq 
\{X_i:i < \sigma\}$.
\end{enumerate}
\mn
As $\alpha \in W(\langle \rangle,\bold S) \Rightarrow \lambda(\langle
\alpha\rangle,S) \le \mu < \lambda$ clearly for some $\theta$
\mn
\begin{enumerate}
\item[$(*)_2$]   $\theta$ is regular $< \mu$
the set $S_1 = \{\delta \in S_0:\lambda(\langle \alpha\rangle,S) =
\theta\}$ is stationary.
\end{enumerate}
\mn
For each $\delta \in S_1$ clearly $\bold x(\langle \delta \rangle)$ is
 a $\theta$-witness and let
\mn
\begin{enumerate}
\item[$(*)_3$]   $(a) \quad X_\delta := \{(\rho,s):\rho \in 
\bold S^{\bold x(<\delta>)}_f$ (i.e. is $\triangleleft$-maximal 
in $\bold S^{\bold x(<\delta>)}$) and $s$

\hskip25pt   is a finite initial segment of
$s^0_{\langle \delta\rangle \char 94 \rho}$ which is a set of members of

\hskip25pt $B_{\langle \delta \rangle}$ of order type $\omega\}$.
\end{enumerate}
\mn
Clearly
\mn
\begin{enumerate}
\item[$(*)_4$]   $X_\delta \subseteq M_\delta$ has cardinality $\theta$.
\end{enumerate}
\mn
We can find $\bar X_\delta$ such that
\mn
\begin{enumerate}
\item[$(*)_5$]   $\bar X_\delta = \langle X_{\delta,i}:i < \sigma\rangle$
is as in $(*)_1$ above, i.e. is $\subseteq$-increasing continuous, so
$X_{\delta,i} \in M_\delta,|X_{\delta,i}| \le \theta$ 
and letting $X_{\delta,\sigma} :=
\cup\{X_{\delta,i}:i < \sigma\}$ we have $X_{\delta,\sigma} = X_\delta$
\sn
\item[$(*)_6$]  $Z_{\delta,i} := \{(\rho,|s|):(\rho,s) \in
X_{\delta,i}\}$ for $i \le \sigma$.
\end{enumerate}
\mn
We define equivalence relation $E$ on $S_1$:
\mn
\begin{enumerate}
\item[$(*)_7$]   $\delta_1 E \delta_2$ iff 
\sn
\begin{enumerate}
\item[$(a)$]  $\bold S^{\langle \delta_1\rangle} = \bold
S^{\langle \delta_2\rangle}$ equivalently
$\bold S_{\bold x(\langle \delta_1\rangle)} = 
\bold S_{\bold x(\langle \delta_2\rangle)}$
\sn
\item[$(b)$]  for each $i < \sigma$ we have $Z_{\delta_1,i}= Z_{\delta_2,i}$
\end{enumerate}
\end{enumerate}
\mn
Clearly $E$ is actually an equivalence relation but $\theta < \mu$ and
 $\mu$ is strong limit hence $E$ has $\le
2^\theta < \mu < \lambda$ equivalence classes.
So for some $\delta^* \in S_1,S_2 := \delta^*/E$ is a stationary
subset of $\lambda$.  Clearly there is $F \in M_0$ which is a one to one
function with range $\subseteq \{\delta:\delta < \lambda,
\delta = \mu \delta$ is $< \lambda$ but $>0\}$ and 
with domain ${\cH}(\lambda)$ hence for every
$\delta \in S_2$ it maps $M_\delta$ onto $\Rang(F) \cap \delta$ 
and \wilog \, if
 $\varrho_1 \triangleleft \varrho_2$ are from ${}^{\lambda >}{\cH}
(\lambda)$ then $F(\varrho_1) < F(\varrho_2)$ and $F(\langle
X_{\delta,j}:j \le i\rangle)$ is $> \sup\{s \cup \Rang(\rho):(\rho,s)
\in X_{\delta,i}\}$.

Let $\alpha^*_{\delta,i} = F(\langle X_{\delta,j}:j \le i\rangle)$ for
$\delta \in S_2,i < \sigma$, so clearly $\langle \alpha^*_{\delta,i}:i
< \sigma\rangle$ is increasing with limit $\delta$ by the choice of
$F$ and of $S_2$, as (see \ref{1n.4}(2) - $\delta = \sup
\cup\{s^{0,\bold x}_\eta:\langle \delta \rangle \trianglelefteq \eta
\in \bold S_f\}$).   By \cite[4.2]{Sh:667}, we can find
$\langle(\nu_{\delta,1},\nu_{\delta,2}):\delta \in S_2\rangle$ such
that
\mn
\begin{enumerate}
\item[$\circledast$]  $(a) \quad \nu_{\delta,\ell} \in {}^\sigma \mu$
\sn
\item[${{}}$]   $(b) \quad \nu_{\delta,1}(i) <
\nu_{\delta,2}(i)$ for $i < \sigma$
\sn
\item[${{}}$]   $(c) \quad$ if $\bold c:\lambda \rightarrow
2^\theta + \sigma$ then for stationarily many $\delta \in S_2$ we have

\hskip25pt  $i < \sigma \wedge \varepsilon < \theta \Rightarrow 
\bold c(\alpha^*_{\delta,i} + \mu \varepsilon + \nu_{\delta,1}(i)) =
\bold c(\alpha^*_{\delta,i} + \mu \varepsilon + \nu_{\delta,2}(i))$.
\end{enumerate}
\mn
Let $\bold y = \bold x \rest \bold S'$ where
$\bold S' = \{\langle \rangle\} \cup \{\rho \in \bold S:\rho \ne
\langle \rangle$ and $\rho(0) \in S_2\}$.

Now at last we shall define the group.  Essentially it will be similar to 
the group $G_{\bold x}$, see \ref{1n.6}(2) so defined inside the
proof of \cite[1.2]{Sh:521} from the system $\bold x$, 
restricted to $\bold S'$ only, but whereas in $(*)^a_{I,\eta}$ before
Fact A in the proof in \cite[1.2]{Sh:521} we use ``$2y_{\eta,m+1} =
y_{\eta,m} + \Sigma\{x[a^\ell_{\eta,m}):\ell < n$ and $a^\ell_{\eta,m}
  \in Y[I]\}$ here we replace $x[a^\ell_{\eta,m}]$ by the difference
  of two, related to $\circledast$; this may become clearer after
  reading the proof.
The $\lambda$-freeness will
be inherited from $\bold S$ being $\lambda$-free.  The non-Whitehead comes from
$\circledast$. 

For $\delta \in S_1$ let $\bold g_\delta:\bold S_{\bold x(\langle
\delta\rangle)} \rightarrow \theta$ be a one-to-one function, so by
the choice of $S_2$ for some $\bold g$ we have $\delta \in S_2
\Rightarrow \bold g_\delta = \bold g$.

As in \cite[\S1]{Sh:521} for each $\eta \in \bold S_f$ and
$\ell < n$ recall, $\langle a^\ell_{\eta,m}:m <
\omega\rangle$ lists $s^\ell_\eta \subseteq B_{\eta \rest (\ell +1)}$, 
let $Y = \cup\{B_\nu:\nu \in \bold S_c\} \backslash 
B^{\bold x}_{\langle \lambda\rangle}$ 
and for $\eta \in \bold S_f,m < \omega$ let $\bold
i_m(\eta)$ be the minimal $i$ such that $(\eta \restriction
[1,n),\{a^0_{\eta,\ell}:\ell \le m\}) \in X_{\eta(0),i}$.  
We define $G$ as the abelian group generated by 

\[
\Xi = \{y_{\eta,m}:m < \omega \text{ and } \eta \in \bold S'_f\} \cup
\{x[a]:a \in Y\} \cup \{z_\beta:\beta < \lambda\}
\]

\mn
freely except for the equations
\mn
\begin{enumerate}
\item[$(*)$]   for $\eta \in \bold S'_f$ and $m < \omega$, so $\delta
 := \eta(0) \in S_2$ the equation (letting $i = \bold i_m(\eta)$)

\begin{equation*}
\begin{array}{clcr}
2y_{\eta,m+1} = y_{\eta,m} &+ \Sigma\{x[a^\ell_{\eta,m}]:0 <\ell < n\} \\
  &+ z_{\alpha^*_{\delta,i} + \mu \bold g(\eta) +
  \nu_{\delta,2}(i)} \\
  & -z_{\alpha^*_{\delta,i} + \mu \bold g(\eta) +
  \nu_{\delta,1}(i)}
\end{array}
\end{equation*}

\mn
recalling $\bold g:\bold S_f \rightarrow \mu$ such that $\bold g
 \restriction \{\eta \in \bold S_f:\eta(0) = \delta\}$ is
one to one, $\alpha^*_{\delta,i} = F(\langle X_{\delta,j}:j \le
 i\rangle)$.  

For $\alpha \le \lambda$ let $G_\alpha$ be the subgroup of $G$
generated by

\begin{equation*}
\begin{array}{clcr}
\{y_{\eta,m}:&m < \omega \text{ and } \eta \in \bold S'_f \text{ and }
\eta(0) < \alpha\} \\
  &\cup\{x[a]:a \in Y \cap \alpha \text{ and } a+1 < \alpha\} \\
  &\cup\{z_\beta:\beta < \alpha \text{ moreover } \beta +1 <
  \alpha\}.
\end{array}
\end{equation*}
\end{enumerate}
\mn
Easily
\mn
\begin{enumerate}
\item[$\oplus_1$]  $G_\alpha$ is a pure subgroup of $G$, increasing
continuous with $\alpha$
\sn
\item[$\oplus_2$]  if $\delta \in S_2$ then $G_{\delta +1}/G_\delta$ 
is isomorphic to $G_{\bold x(\langle \delta\rangle)}$
which is not Whitehead.
\end{enumerate}
\mn
[Why?  By clause (d) of the assumption; now at first glance the set of
  generators and equations in the proof of \cite[1.2]{Sh:521} and in
  this proof are different.  But note that only $y_\eta,\langle \alpha
  \rangle \triangleleft \eta \in \bold S_f$ and $x[a],a \in
  \cup\{s^\ell_\eta:\ell \in [1,n),\langle \alpha \rangle
    \trianglelefteq \eta \in \bold S_f\}$ appear in the equation.
    Alternatively use Remark \ref{a8} and prove $G_{\alpha
      +1}/G_\alpha$ is isomorphic to $G_{\bold x,\alpha +1}/G_{\bold
      x,\alpha}$ in \cite[1.2]{Sh:521} notation; again note that the
    $z_\alpha$ - here and $x[a],a \in \cup\{s^{0,\bold x}_\eta:\eta
    \in \bold S_f\}$ disappear.]

But recall
\mn
\begin{enumerate}
\item[$\oplus_3$]   $G$ is strongly $\lambda$-free, moreover if
$\alpha \in \lambda \backslash S_2$ and $\beta \in (\alpha,\lambda)$
 then $G_\beta/G_\alpha$ is free.
\end{enumerate}
\mn
[Why?  As in the proof of Fact A inside the proof of
\cite[1.2]{Sh:521}.]
\mn
\begin{enumerate}
\item[$\oplus_4$]  $G$ is not Whitehead.
\end{enumerate}
\mn
[Why?  We choose $(H_\alpha,h_\alpha,g_\alpha)$ by induction on $\alpha \le
\lambda$ such that
\mn
\begin{enumerate}
\item[$(a)$]  $H_\alpha$ is an abelian group extending $\bbZ$
\sn
\item[$(b)$]  $h_\alpha$ is a homomorphism from $H_\alpha$ onto
$G_\alpha$ with kernel $\bbZ$
\sn
\item[$(c)$]  $g_\alpha$ is a function from $G_\alpha$ to
$H_\alpha$ inverting $h_\alpha$ (but in general not a homomorphism)
\sn
\item[$(d)$]  $H_\alpha$ is increasing continuous with $\alpha$
\sn
\item[$(e)$]  $h_\alpha$ is increasing continuous with $\alpha$
\sn
\item[$(f)$]  $g_\alpha$ is increasing continuous with $\alpha$
\sn
\item[$(g)$]  if $\alpha = \delta +1$ and $\delta \in S_2$ then
there is no homomorphism $g^*$ from $G_\alpha$ into $H_\alpha$
inverting $h_\alpha$ such that:
\begin{enumerate}
\item[$\odot$]   $i < \sigma \wedge \varepsilon < \theta
\Rightarrow g^*(z_{\alpha^*_{\delta,i} + \mu \varepsilon +
\nu_{\delta,1}(i)}) - g_{\delta +1}(z_{\alpha^*_{\delta,i} + \mu
\varepsilon + \nu_{\delta,1}}) = g^*(z_{\alpha^*_{\delta,i} + \mu
\varepsilon + \nu_{\delta,2}(i)}) - g_{\delta +1}
(z_{\alpha^*_{\delta,i} + \mu \varepsilon + \nu_{\delta,2}(i)})$
\end{enumerate}
\end{enumerate}
\mn
(note: the subtraction in $\bbZ$, the kernel of $h_\alpha$)

For $\alpha = 0,\alpha$ limit and $\alpha = \beta +1,\beta \notin S_2$
this is obvious.  For $\alpha = \delta +1,\delta \in S_2$ it is known
that if instead of $\odot$ in clause (g) we know $g^* \rest G_\delta$ this
is possible.  But $\odot$ gives all the necessary information.  In
more details let $G'_\delta$ be the subgroup of $G_\delta$ generated
by $\{z_{\alpha^*_{\delta,i} + \mu \varepsilon + \nu_{\delta,2}(i)} -
z_{\alpha^*_{\delta,i} + \mu \varepsilon +
\nu_{\delta,1}(i)}:\varepsilon < \theta$ and $i < \sigma\}$.

Let $G'_{\delta +1}$ be the subgroup of $G_{\delta +1}$ generated by
$G'_\delta \cup \{y_{\eta,m}:\langle \delta \rangle \trianglelefteq
\eta \in \bold S'_f\}$.  Clearly $G'_\delta$ is a pure subgroup of
$G'_{\delta +1}$ and of $G_\delta$ and $G_{\delta +1} = G'_{\delta +1}
{\underset{G'_\delta} \oplus} G_\delta$.

Let $H'_\delta = h^{-1}_\delta(G'_\delta)$, clearly $h_\delta \rest
H'_\delta$ is a homomorphism from $H'_\delta$ onto $G'_\delta$ with
kernel $\bbZ$.

Clearly 
\mn
\begin{enumerate}
\item[$\oplus_{4.1}$]  if $g',g'' \in \text{ Hom}(G_\delta,H_{\delta +1})$
invert $h_\delta$ and both satisfies $\odot$ of clause (g) then
$g' \rest G'_\delta = g'' \rest G'_\delta$.
\end{enumerate}
\mn
So $|{\cG}_\delta| \le 1$ where

\begin{equation*}
\begin{array}{clcr}
{\cG}_\delta = \{g \rest G'_\delta:&g \text{ is a homomorphism
from } G_\delta \text{ to } H_\delta \\
  &\text{ inverting } h_\alpha \text{ and satisfying } \odot 
\text{ in clause (g)}\}.
\end{array}
\end{equation*}

\mn
Let $g^*$ be the unique member of $\cG$ if $\cG$ is non-empty and
otherwise let it be any homomorphism from $G'_\delta$ into $H'_\delta$
inverting $h_\delta \rest H'_\delta$, exist as $G'_\delta$ is free.

We now choose $(H'_{\delta +1},h'_{\delta +1})$ such that
\mn
\begin{enumerate}
\item[$\bullet$]  $H'_\delta \subseteq H'_{\delta +1}$
\sn
\item[$\bullet$]  $h'_{\delta +1} \in \Hom(H'_{\delta +1},G'_{\delta +1})$
\sn
\item[$\bullet$]  $h'_{\delta +1}$ has kernel $\bbZ$
\sn
\item[$\bullet$]  $h'_{\delta +1}$ extends $h_\delta \rest G'_\delta$
\sn
\item[$\bullet$]  $g^*$ cannot be extended.
\end{enumerate}
\mn
Next \wilog \, $H'_{\delta +1} \cap H_\delta = H'_\delta$, let
$H_{\delta +1} = H'_{\delta +1} \underset{H'_\delta} \oplus
H_\delta$, and let $h_{\delta +1} \in \Hom(H_{\delta +1},G_{\delta
  +1})$ extend $h_\delta,h'_{\delta +1}$.  Let $g_{\delta +1}
\supseteq g_\delta$ invert $h_{\delta +1}$ but is not necessarily a
homomorphism. So we have chosen
$(H_\alpha,h_\alpha,g_\alpha)$ such that: if 
${\cG}_\delta \ne \emptyset$ then we cannot find $g' \in 
\Hom(G'_{\delta +1},\bbZ)$ inverting $h_\alpha$ such that $g'
\rest G'_\delta \in {\cG}_\delta$.  This suffices for carrying out the
induction.

Having carried the induction, clearly $h = h_\lambda$ is a
homomorphism from $H = H_\lambda$ onto $G_\lambda = G$ with kernel
$\bbZ$.  To show that $G$ is not a Whitehead group it suffices to
prove that $h$ is not invertible as a homomorphism.  But if $g \in
\text{ Hom}(G,H)$ inverts $h$ then $x \in G \Rightarrow
g(x)-g_\lambda(x) \in \bbZ$.

We define a function $\bold f$ with domain $\lambda$: for $\alpha <
\lambda,\zeta < \mu:\bold f(\mu \alpha + \zeta) = \langle g(z_{\mu^2
\alpha + \mu \varepsilon + \zeta}) - g_\lambda(z_{\mu^2 \alpha +
\mu \varepsilon + \zeta}):\varepsilon < \theta\rangle$.

So for stationarily many $\delta \in S_2$ we have $i < \sigma
\Rightarrow \bold f(\alpha^*_{\delta,i} + \nu_{\delta,1}(i)) = \bold
f(\alpha^*_{\delta,i} +\nu_{\delta,2}(i))$.   For any such $\delta$ we
get a contradiction by clause (g) of the construction, so we have proved
$\oplus_4$.

This finishes the proof of part (1), as $G = G_\lambda$ is as
required.

\noindent
2) For the proof of part (2) can use:
\mn
\begin{enumerate}
\item[$\boxtimes$]  for regular uncountable $\lambda$, the following
conditions are equivalent
\sn
\begin{enumerate}
\item[$(a)$]  every $\lambda$-free abelian group of cardinality
$\lambda$ is Whitehead
\sn
\item[$(b)$]  for every $\lambda$-free abelian group of
cardinality $\lambda$ we have Hom$(G,\bbZ) \ne 0$.
\end{enumerate}
\end{enumerate}
\mn
[Why?  If (a) and $G$ as in(b), let $h$ be a pure embedding of $\bbZ$
  into $G$, let $G' = G/\text{Rang}(h)$ and use the definition of
  ``$G'$ is a Whitehead group".  If $G,H$ and $h \in \Hom(H,G)$ form a
  counterexample then we can find a purely increasing continuous
  sequence $\langle G_\alpha:\alpha \le \lambda\rangle$ and $\langle
  x_\alpha:\alpha < \lambda\rangle$ such that: $\{x_\alpha:\alpha <
  \lambda\} = \{x \in G_\lambda:\bbZ x$ is a pure subgroup of
  $G_\lambda\}$ and for each $\alpha$, there is a pure embedding
  $h_\alpha$ of $H$ into $G_{\alpha +1}$ such that $h(1_{\bbZ}) =
  x_\alpha$ and $G_{\alpha +1} = G_\alpha \underset{\bbZ x_\alpha}
  \oplus h_\alpha(H)$.

Easily $G_\lambda$ contradicts clause (b).]

We can also construct directly.
\end{PROOF}
\newpage

%\bibliographystyle{alphacolon}
%\bibliography{lista,listb,listx,listf,liste,listz}

\end{document}